\documentclass[12pt,leqno]{amsart}
\usepackage{amssymb}
\usepackage{amscd}
\usepackage{amscd}
\usepackage{amsmath}
\usepackage{amsthm}

\newcommand{\Q}{\mathbb Q}
\newcommand{\Z}{\mathbb Z}
\newcommand{\C}{\mathbb C}
\newcommand{\R}{\mathbb R}
\newcommand{\A}{\mathbb A}

\newcommand{\F}{\mathbb F}
\newcommand{\G}{\mathfrak G}
\newcommand{\hH}{\mathcal H}
\newcommand{\K}{\mathcal K}
\newcommand{\fO}{\mathcal O}
\DeclareMathOperator{\Todd}{Todd}
\DeclareMathOperator{\td}{td}
\def\id{{\sl id}}

\makeatletter
\def\1@section{\@tocline{1}{4pt}{1pc}{}{}}
\def\1@subsection{\@tocline{2}{0pt}{2pc}{5pc}{}}
\makeatother

\begin{document}
\title{\bf MODULAR FORMS AND CALABI-YAU VARIETIES}
\author[Kapil Paranjape and Dinakar Ramakrishnan]{Kapil Paranjape$^{1}$ and Dinakar Ramakrishnan$^{2}$}
\thanks{$^1$ Partly supported by the the JCBose Fellowship
 of the DST (IISERM:10-DST-JCB-F.3) and IMSc, Chennai}
\thanks{$^2$ Partly supported by the NSF grants DMS-0701089 and DMS-1001916}
\address{Department of Mathematical Sciences \\ IISER\\ Mohali, India}
\email{kapil.iiserm@gmail.com}
\address{Department of Mathematicss \\ California Institute of Technology \\
Pasadena, CA 91125} \email{dinakar@caltech.edu}

\date{}

\setcounter{page}{1} \maketitle

\bigskip

\section*{Introduction}

\medskip

Let $f(z)=\sum\limits_{n=1}^\infty a_nq^n$ be a holomorphic newform
of weight $k \geq 2$ relative to $\Gamma_1(N)$ acting on the upper
half plane $\hH$. Suppose the coefficients $a_n$ are all rational.
When $k=2$, a celebrated theorem of Shimura asserts that there corresponds an elliptic
curve $E$ over $\Q$ such that for all primes $p\nmid N$, $a_p =
p+1-\vert E(\F_p)\vert$. Equivalently, there is, for every prime
$\ell$, an $\ell$-adic representation $\rho_\ell$ of the
absolute Galois group $\G_\Q$ of $\Q$, given by its action on the
$\ell$-adic Tate module of $E$, such that $a_p$ is, for any $p\nmid
\ell N$, the trace of the Frobenius $Fr_p$ at $p$ on $\rho_\ell$.

The primary aim of this article is to provide
some positive evidence for the expectation of Mazur and
van Straten that for every $k\geq2$, any (normalized) newform
$f$ of weight $k$ and level $N$ should, if it has rational coefficients,
have an associated Calabi-Yau
variety $X/\Q$ of dimension $k-1$ such that
\begin{enumerate}
\item[(Ai)]
The $\{(k-1,0),(0,k-1)\}$-piece of $H^{k-1}(X)$ splits off as a
submotive $M_f$ over $\Q$,
\item[(Aii)] $a_p={\rm tr}(Fr_p\, \vert \,
M_{f,\ell})$, for almost all $p$, and
\item[(Aiii)]det$(M_{f,\ell})=\chi_\ell^{k-1}$,
\end{enumerate}
where $\chi_\ell$ is the $\ell$-adic cyclotomic character.

\medskip

In fact we expect (Aiii) to hold for every $p$ not dividing $\ell N$.

\medskip

\noindent{}Furthermore, we even hope that in addition the following holds:

\smallskip

\noindent \, \, (Aiv) $X$ admits an involution $\tau$ which acts by $-1$ on
$H^0(X,\Omega^{k-1})$.

\smallskip

We anticipate that the involution $\tau$ can be chosen in such a way to make the quotient
$X/\tau$ a rational variety.

This typically holds in our examples below. This extra structure
is natural to want and is needed for understanding a variety of operations like products
and twists.
Such a $\tau$ obviously exists for $k=2$, in which case $X$ is (by Shimura) an
elliptic curve over $\Q$, given by an equation $y^2=f(x)$, and
the involution $\tau$ sends $(x,y)$
to $(-x,y)$, thus acting by $-1$ on the holomorphic differential
$\omega = dx/y$ which spans $H^0(X,\Omega^1)$. (Of course $X/\tau$ is in this case ${\mathbb P}^1$.)

To fix ideas,
one could think of a motive as a semisimple motive $M$
relative to absolute Hodge cycles (\cite{DMOS}), with avatars
$(M_B,M_{\rm dR},M_\ell)$ of Betti, de Rham and $\ell$-adic realizations. Since every Calabi-Yau manifold of dimension $1$ is an elliptic
curve, (Ai) through (Aiv) provide a natural extension of what one has for $k=2$. Of
course for $k
>2$, one knows by Deligne (\cite{De}) that there is an irreducible,
$2$-dimensional $\ell$-adic representation $\rho_\ell$ of Gal$(\overline \Q/\Q)$
so that
$(Aii), (Aiii)$ hold with $M_{f,\ell}$ replaced by $\rho_\ell$.

In the first part of this article we will focus on the forms $f$ of
even weight case for small levels, before moving on in the second part to formulate an analogue
for regular selfdual cusp forms on GL$(n)/\Q$ and see how
the framework is compatible with the principle of functoriality.

When the weight is odd, the
$\Q$-rationality forces $f$ to be of CM type, and for weight $3$, we refer to
the paper of Elkies and Schuett (\cite{ES}) for a beautiful result.

For non-CM newforms $f$ of weight $k>2$ with $\Q$-coefficients, we in
fact hope for more, namely that the cohomology ring of $X$ will be spanned by $M$ and
the Hodge/Tate classes of various degrees; in particular, $X$ should
be rigid in this case.

It is a difficult problem in dimensions $>3$ to find smooth models
of varieties with trivial canonical bundles, and for this reason we
formulate our questions for such varieties with mild singularities.
By a {\it Calabi-Yau variety} over a field $k$, we will mean an
$n$-dimensional projective variety $X/k$ on which the canonical
bundle $\K_X$ is defined such that
\begin{enumerate}
\item[(CY1)] $\K_X$ is trivial; \, and
\item[(CY2)] $H^m(X, \fO_X) = 0$ for all (strictly) positive $m <n$.
\end{enumerate}
More precisely, we will want such an $X$ to be normal and Cohen-Macaulay,
so that the dualizing sheaf $\K_X$ is defined, with the singular
locus in codimension at least $2$, so that $\K_X$ defines a Weil
divisor; finally, $X$ should be Gorenstein, so that $\K_X$ will
represent a Cartier divisor. Ideally we would like to singular locus
$X_{\rm sing}$ to be of dimension $\leq \left[\frac{n-1}{2}\right]$.

In addition to these properties, we would ideally also like $X$ to be realized as a double cover
$\pi: X \to Y$, where $Y$ is a projective smooth rational variety with negative ample canonical divisor.
(Such an $X$ will be automatically Gorenstein so that $\K_X$ is defined, and
$\pi_\ast(\fO_X)=\fO_Y\oplus L$,
for a line bundle $L$ on $Y$ with $L^2$ giving the branch locus; moreover,
$\pi_\ast(\K_X) = \K_Y \oplus \K_Y\otimes L^{-1}$, forcing $L=\K_Y$.)

\medskip

Here is our first result:

\medskip

\noindent{\bf Theorem 1} \, \it Fix $\Gamma = \Gamma_1(N)$, $N \leq
5$. Let $k$ be the first even weight s.t. dim$(S_k(\Gamma)) =1$.
Then $\exists$ a Calabi-Yau variety $V(f)/\Q$ with an involution $\tau$
associated to the new
generator $f$ of $S_k(\Gamma)$ satisfying $(Ai)$ and $(Aiv)$. In
fact, when $N\leq 5$, $V = V(f)$ is birational over $\Q$ to the Kuga-Sato
variety $\tilde{\mathcal E}^{(k-2)}_N$.  \rm

\medskip

In particular, this applies to the Delta function $\Delta(z)
= \sum_{n=1}^\infty \tau(n)q^n$ $= \prod_{m\geq 1}(1-q^m)^{24}$, for
$N=1, k=12$. By the Kuga-Sato variety, we mean a smooth
compactification (over $\Q$) of the fibre product $E^{k-2}_N$ of
the universal elliptic curve $E_N$ over the model over $\Q$ of the
modular curve $\Gamma_1(N)\backslash\hH$.

\medskip

It is well known that $S_k(\Gamma_1(N))$ is {\bf one-dimensional}
when $(N,k)$ equals ${\bf(1,12)}, (1,16), {\bf (2,8)}, (2,10)$,
${\bf (3,6)}$, $(3,8)$, ${\bf (4,6)}$, ${\bf (5,4)}$, $(5,6)$, ${\bf
(6,4)}, (7,4)$. (The ones in bold are the cases to which the Theorem
applies.) For example, for the case ${\bf(2,8)}$, the generator is
$$
f(z) \, = \, q-8q^2+12q^3-210q^4+1016q^5 +\dots
$$

\medskip

Recall that a newform $f(z)=\sum_{n=1}^\infty a_nq^n$ is of {\it
CM-type} iff there is an odd, quadratic Dirichlet character $\delta$
such that $a_p=a_p\delta(p)$ for almost all primes $p$.
Equivalently, if $K$ is the imaginary quadratic field cut out by
$\delta$, $a_p=0$ for all $p$ which are inert in $K$.

As a first step, we may ask for a {\it potential statement}, i.e., the
association, over a finite extension $F$ of $\Q$, of a Calabi-Yau variety $V/F$
to a newform $f$ with $\Q$-coefficients. Here we state
a modest result in this direction,
already known in different ways, just to show that it fits into our framework:

\medskip

\noindent{\bf Proposition 2} \, Let $f$ be a newform of weight $k \geq
3$ of CM type with rational coefficients. Then $\exists$ a
Calabi-Yau $(k-1)$-fold $X$ defined over a number field $F$ such
that (Ai), (Aii) and (Aiii) hold over $F$. This $X$ arises as a Kummer
variety associated to an elliptic curve $E$ with complex
multiplication. When $k\leq 4$, $X$ can be taken to be a smooth model.\rm

\medskip

For $k=4$, all but one form $f$ have been treated in Cynk-Schuett (\cite{Cy-Sch}).

\medskip

Again, when $k=3$, there is a much more precise and satisfactory result
over $\Q$ in the work of Elkies and Schuett \cite{ES}.

\medskip

In the converse direction, if $M$ is a simple motive over $\Q$ of
rank $2$, with coefficients in $\Q$, of Hodge type $\{(w,0),(0,w)\}$
with $w>0$, then the general philosophy of Langlands, and also a
conjecture of Serre, predicts that $M$ should be modular and be
associated to a newform $f$ of weight $w+1$ with rational
coefficients. This is part of a very general phenomenon, and applies to
motives occurring in the cohomology of smooth projective varieties over $\Q$.
In any case, it applies in
particular to Calabi-Yau threefolds over $\Q$ whose
$\{(3,0),(0,3)\}$-part splits off as a submotive. In this context,
there have been a number of beautiful results, some of which have been
described in the monographs \cite{CYMS} and \cite{MSV}. See also \cite{M},
\cite{GKY} and
\cite{G-Y}. They are entirely consistent
with what we are trying to do in the opposite direction, and also
provide supporting examples. It should perhaps be remarked that these results
(relating to the modularity of rigid Calabi-Yau threefolds) can now be deduced
{\it en masse} from the proof of Serre's conjecture due
to Khare and Wintenberger (\cite{KhW}), with a key input from Kisin (\cite{Ki}).
See also the article of Dieulefait (\cite{D}).

\medskip

Let us now move to a more general situation. Fix any positive
integer $n$ and suppose that $f$ is a (new) Hecke eigen-cuspform on
the symmetric space
$$
{\mathcal D}_n: = {\rm SL}(n, \R)/{\rm SO}(n),
$$
relative to a congruence subgroup $\Gamma$ of SL$(n, \Z)$, which is
{\it algebraic} and {\it regular}. For $n=2$, $f$ is algebraic and
regular iff it is holomorphic of weight $\geq 2$. In general, one
considers the cuspidal automorphic representation $\pi$ of
GL$(n,\A)$ which is generated by $f$, and by Langlands the
archimedean component $\pi_\infty$ corresponds to an $n$-dimensional
representation $\sigma_\infty$ of the real Weil group $W_\R$, which
contains $\C^\ast$ as a subgroup of index $2$. One says (\cite{C})
that $\pi$ is algebraic if the restriction of $\sigma_\infty$ to
$\C^\ast$ is a sum of characters $\chi_j$ of the form $z\to
z^{p_j}\overline z^{q_j}$, with $p_j, q_j \in \Z$, and it is regular
iff $\chi_i\ne\chi_j$ when $i\ne j$. Such an $f$ contributes to the
{\it cuspidal cohomology} $H^\ast_{\rm cusp}(\Gamma\backslash
{\mathcal D}_n, V)$ relative to a {\it local coefficient system} $V$
in a specific degree $w=w(f)$. Moreover, $f$ is rational over a
number field $\Q(f)$, defined by the Hecke action on cohomology,
which preserves the cuspidal part ({\it loc. cit.}). There is
conjecturally a motive $M(f)$ over $\Q$ of rank $n$, with
coefficients in $\Q(f)$, and weight $w$. By Clozel \cite{C2}, $M_\ell(f)$
exists for suitable $f$ which are in addition {\it essentially
selfdual}.

Now let $f$ be an algebraic, regular, essentially selfdual newform
of weight $w$ relative to $\Gamma\subset {\rm SL}(n,\Z)$, with
$L$-function $L(s,f)=\prod_p L_p(s,f)$, such that $\Q(f) = \Q$. Then
our question is if there exists a Calabi-Yau variety $X/\Q$ of
dimension $w$ with an involution $\tau$ such that
\begin{enumerate} \item[(Ci)] There is a submotive $M(f)$ of
$H^w(X)$ of rank $n$ such that $M(f)^{(w,0)}=H^{w,0}(X)$,
\item[(Cii)] $L_p(s,f)=L_p(s,M_\ell(f))$ for almost all $p$, \, and
\item[(Ciii)]the quotient of $X$ bt $\tau$ is a rational variety.
\end{enumerate}
where $L_p(s,M_\ell(f))$ equals, at any prime $p\ne \ell$ where the
$\ell$-adic realization $M_\ell(f)$ is unramified, det$(I-Fr_pp^{-s}
\, \vert \, M_\ell(f))^{-1}$. Again we would like to be able to find
an $X$ having good reduction outside the primes dividing $tN$, where
$n$ is the level of $f$ and $t$ the order of torsion in $\Gamma$,
such that $(Cii)$ holds for any such $p \ne \ell$.

Thanks to the {\it principle of functoriality}, one should be able
to obtain a certain class of $\Q$-rational, regular, algebraic,
essentially selfdual newforms $f$ by transferring forms on (the
symmetric domains of) smaller reductive $\Q$-subgroups $G$ of
GL$(n)$. The simplest instance of this phenomenon is given by the
symmetric powers sym$^m(g)$ of classical $\Q$-rational, non-CM
newforms $g$ of weight $k$. One knows by Kim and Shahidi (\cite{KSh1}, \cite{K})
that for $m\leq 4$, $f={\rm sym}^m(g)$ is a cusp form on GL$(m+1)$.
(Recently, this has been extended to $m=5$ (and further) in  the works of Clozel and Thorne,
and of Dieulefait.)
Here is our third result:

\medskip

\noindent{\bf Theorem 3} \, \it Let $g$ be a non-CM, elliptic
modular newform of weight $2$, level $N$ and trivial character,
whose coefficients $a_n$ lie in $\Q$. Then for for any $m >0$, there
is a Calabi-Yau variety $X_m$ with an involution $\tau$ over $\Q$ of dimension $m$ associated
to $(g, {\rm sym}^m)$ such that $(Ci), (Cii), (Ciii)$ hold relative to
$M_\ell={\rm sym}^m(\rho_\ell(g))$. Moreover, for $m\leq 3$,
$X_m$ can be taken to be non-singular, with good reduction outside
$N$. \rm

\medskip

Here $X_2$ is just the familiar Kummer surface attached to $E\times
E$, where $E=X_1$ is the elliptic curve$/\Q$ defined by $g$. But the
case $m=3$ is interesting, especially since it is not rigid, thanks
to the Hodge type being $\{(3,0),(2,1),(1,2),(0,3)\}$, with each
Hodge piece being one-dimensional. In fact, in that case, sym$^3(g)$
corresponds (by \cite{RaSh}) to a holomorphic Siegel modular cusp form
$F$ of genus $2$ and (Siegel) weight $3$. Such an $F$ contributes to
the cohomology in degree $3$ of the Siegel modular threefold $V$ of
level $N^3$. Since the geometric genus of $V$ is typically $>1$, it
cannot be Calabi-Yau. However, there should be, as predicted by the
Hodge and Tate conjectures, an algebraic correspondence between
$X_3$ and $V$ (for any $N$).

\medskip

One also knows (cf. \cite{Ra2}) that given two non-CM newforms $g, h$ of
weights $k, r \geq 2$ respectively, then there is an algebraic
automorphic form $f=g\boxtimes h$ on GL$(4)/\Q$, which will be
cuspidal and regular if $k\ne r$. If $g, h$ are $\Q$-rational, then
so is $f$. Moreover, $f$ is essentially selfdual because $g$ and $h$
are.

\medskip

\noindent{\bf Theorem 4} \, \it Let $g,h$ be $\Q$-rational, non-CM
newforms as above of respective weights $k, r >1$, with $k\ne r$.
Suppose we have Calabi-Yau varieties with involutions $(X(g), \tau_g),
(X(h), \tau_h)$ over $\Q$ attached to $g,
h$ respectively, satisfying $(Ai)$ through $A(iv)$. Put $f=g\boxtimes h$, so
that $w(f)=(k-1)(r-1)$. Then there is a Calabi-Yau variety with involution $(X(f), \tau_f)$
over $\Q$
of dimension $w(f)$ such that $(Ci)$ through $(Ciii)$ hold. \rm

\medskip

The point is that the product $Z:=X(g)\times X(h)$ has the desired
submotive in degree $w$, but it has global holomorphic $m$-forms for
$m=k-1$ and $m=r-1$. We exhibit an involution $\tau$ on $Z$ such
that when we take the quotient by $\tau$, these forms get killed and
we get a Calabi-Yau variety with reasonable singularities. To get
unconditional examples of this Theorem, take $k=2$ and choose $h$ to
be one of the examples of Theorem $A$ of weight $r>2$. To be
specific, we may take $g$ to be the newform of weight $2$ and level
$11$, and $h$ to be the newform of weight $4$ and level $5$, in
which case $f=g\boxtimes h$ has level $55^2$, and $X(f)$ is a
Calabi-Yau fourfold.
A similar inductive construction of Calabi-Yau varieties is also found in a paper of Cynk and Hulek
(\cite{CyH}).

\medskip

In sum, it is a natural question, given Shimura's work on forms of weight $2$,
if there are Calabi-Yau varieties with an involution associated to forms of higher weight
with rational coefficients, and a preliminary version of this circle of questions was
raised by the first author
in a talk at the Borel memorial conference at Zhejiang University in
Hangzhou, China, in 2004, and quite appropriately, Dick Gross, who was in
the audience, cautioned against hoping for too much without sufficient evidence.
Over the past years, there has been some positive evidence, though small, and
even if there is no $V$ in general, especially for non-CM forms $f$ of even weight,
the examples where one has nice Calabi-Yau varieties $V$ enriches them considerably, and
one of our aims is to understand the Hecke eigenvalues $a_p$ in such cases a bit better in terms of
counting points of $V$ mod $p$. Since then we have learnt from \cite{ES} (and Noriko Yui) that
the question of existence of a Calabi-Yau variety $V$ associated to $f$ were earlier raised by
Mazur and van Stratten. What we truly hope for is that in addition,
$V$ will be equipped with an involution $\tau$ acting by $-1$ on the unique global holomorphic form of maximal degree, so one can form products, etc., and also deal with quadratic twists.
We have some interesting examples in the $3$-dimensional case (where $k=4$), and a sequel to this paper will also contain a discussion of these matters, as well as a way to get nicer models in certain higher
dimensional examples.

\medskip

The Ramanujan coefficients $\tau(p)$ of the Delta function have been a source of much research. Our own work was originally motivated by the desire to express them in terms of the zeta function of a Calabi-Yau variety, though our path has diverged somewhat.  In a different direction, a very interesting monograph of Edixhoven, et al (\cite{EC}) yields a deterministic algorithm for computing $\tau(p)$ with expected running time which is polynomial in $\log p$. They do this by relating the associated Galois representation mod $\ell$ to the geometry of certain effective divisors on the modular curve $X_1(5\ell)$.

\medskip

We thank all the people who have shown interest in this work, and a special thanks must go to Matthias Schuett who made a number of useful comments on the earlier version (2008) as well as a recent version (of two months ago). One of us (K.P.) would like to thank the JCBose Fellowship of the DST (IISERM:10-DST-JCB-F.3) as well as Caltech and the IMSc, Chennai, where some of the work was done, and the second author (D.R.) would like to thank the NSF for continued support through
the grants DMS-0701089 and DMS-1001916, and also IMSc, Chennai, for having him visit at various times.

\medskip

\section{An intuitive picture of the geometric construction}

\medskip

In this section we will give an idea behind our construction of varieties with trivial
canonical bundles arising as birational models of elliptic modular
varieties $V$ with non-positive canonical bundles. (S.T.~Yau has informed us that he earlier had a similar construction, albeit in a different context.) We will use the intuitive language
of divisors and linear systems. The content here will not be used
in the succeeding sections, where we will use sheaves and do everything precisely in the
different cases at hand.

Recall that $V$
arise as fibre products of the universal family of elliptic curves
$E$ with additional structures over the modular curve associated to
a congruence subgroup $\Gamma$ of SL$(2,\Z)$. When $\Gamma$ is
$\Gamma_1(N)$ (resp. $\Gamma_0(N)$), the additional structure is a
point (resp. subgroup) of order $N$. We have a slight preference for
$\Gamma_1(N)$ over $\Gamma_0(N)$, because in the latter case, one
gets, due to the existence of $-I$, only a coarse moduli space.

\medskip

Suppose we want to parametrize triples of the form $(E,S;R)$
where $E$ is a curve of genus 1, $S$ a finite set of points on $E$,
and $R$ a finite set of linear equivalence relations on the points
$S$. We will think of $E$ as a curve of genus 1 without a specified origin
(by making an identification with line bundles of degree 1 on it),
hence with no group structure, and the relations in $R$ are taken to
hold in the divisor class group. (Once we pick a point $o$ of $E$, we can
of course identify this class group with $E \times \Z$ by sending $(e,n)$ to
$\left((e-o)+n.o\right)$.) We will take at least a portion of $S$ to consist of general points,
and we assume that $R$ contains all the relations between the chosen
points. Additional relations may hold
for special triples (E,S,R) but generically, only those in $R$
will hold.

Suppose also that there is a surface $X$, and a linear system $P$ of
divisors linear equivalent to $-K_X$, where $K_X$ is the canonical
divisor of $X$. Assume that for a ``general'' datum $(E,S;R)$ we
have a uniquely determined element of $P$, so that the projective
space $P$ parametrizes triples as above up to birational isomorphism.
Here by ``general datum'' we mean a point in an open subset of the
parameter space (or moduli stack to be precise).

Let $W\subset\Gamma(X,\fO(-K_X))$ be the subspace so that $P$ is the
associated projective space and $n={\rm rank} W$. We have a natural
homomorphism
\[   \phi: V\otimes \fO_{X^n} \to \oplus_{i=1}^n \, pr_i^* \fO(-K_X), \]
where $pr_i$ denotes the projection onto the $i$-th factor.
The divisor $D$ where the determinant $\det(\phi)$ vanishes parametrizes,
birationally, tuples of the form
$(E,S;R; p_1,\dots,p_n)$ where the $p_i$ are $n$ additional points which
are not subject to any additional relation. Moreover, $D$ has
trivial canonical bundle.

If $n$ is at most the dimension of the linear system, then the parameter space is rational.
But new things happen when $n$ is larger than that, when one gets a divisor on a product of rational surfaces, in fact given by the vanishing of det$(\phi)$. This moduli space $V_n$, say, fibers over the rational variety $V_{n-1}$, with the general fibre being an elliptic curve. The involution $x \mapsto -x$ on the general fibre gives rise to an involution $\tau$ on $V_n$. So $V_n/\tau$ fibers over a rational variety with ${\mathbb P}^1$ fibres. Hence $V_n/\tau$ is unirational.

\medskip

\section{The modular varieties of interest}

In the context of the congruence subgroup $\Gamma_1(N)$ of
SL$(2,\Z)$, the elliptic modular varieties $V$ with non-positive canonical bundles
are associated with pairs $(N,k)$ such that there is at most one
modular form of level $N$ and weight $k$. The complete list of such
pairs is given in the following table:

\begin{center}
\begin{tabular}{|c|c|}
\hline
$N$ & $k$\\
\hline
1   & $\leq 23$, 25, 26, odd\\
\hline
2   & $\leq 11$, odd\\
\hline
3   & $\leq 8$\\
\hline
4   & $\leq 6$\\
\hline
5, 6   & $\leq 4$\\
\hline=
7, 8   & $\leq 3$\\
\hline
9, 10, 11, 12, 14, 15   & 2\\
\hline
\end{tabular}
\end{center}

\medskip

One can similarly make the corresponding tables for the groups
$\Gamma_0(N)$ and $\Gamma_0(N^2)\cap\Gamma_1(N)$.

\bigskip

\section{The Calabi-Yau $11$-fold associated to $\Delta$}

\medskip

$\Delta$ is a generator of $S_{12}({\rm SL}(2, \Z))$. The object is
to show that the Kuga-Sato variety ${\mathcal E}^{(10)}$ is
birational to an eleven-dimensional Calabi-Yau variety $V$. We will
use $\equiv$ to denote birational equivalence. It is easy to see
that for any $r \geq 0$,
$$
{\mathcal E}^{(r)} \, \equiv \, {\mathcal M}_1(r+1),
$$
where ${\mathcal M}_g(k)$ is the {\it moduli space of genus $g$ curves
with $k$ marked points}, with compactification $\overline{\mathcal
M}_g(k)$.

\medskip

So we need to find a {\it birational model} $V$ of
$\overline{\mathcal M}_1(11)$ such that $V$ is Calabi-Yau. Let $S$
be the surface obtained by {\it blowing up $4$ general points} $P_1,
P_2, P_3, P_4$ in ${\mathbb P}^2$. Let $E$ be an elliptic curve with
$n+5$ general points $Q_0, Q_1, \dots, Q_{n+4}$, and use $\vert
3Q_0\vert$ to define a morphism $E \to {\mathbb P}^2$.

\medskip

Using an automorphism of ${\mathbb P}^2$ we may assume: $Q_i = P_i$
for $1\leq i\leq 4$. The embedding $E \to {\mathbb P}^2$ lifts to a
morphism $\varphi: E \to S$, and the {\it adjunction formula} gives
$\varphi(E) \in \vert {\mathcal K}_S^{-1}\vert$.

\medskip

Get a {\it rational map} ${\mathcal M}_1(n+5) \rightarrow S^n,$
$$(E, \{Q_0, \dots,
Q_{n+4}\}) \to (P_5, \dots, P_{n+4}).
$$
$W:= \Gamma(S, {\mathcal K}_S^{-1})$ has dim. $6$, and $\exists$ a
hom (of sheaves on $S^n$):
$$
f_n: \, W \otimes {\mathcal O}_{S^n} \, \rightarrow \, {\mathcal
K}_S^{-1} \boxtimes \dots \boxtimes {\mathcal K}_S^{-1}.
$$

\medskip

Ker$(f_n)$ is the vector space of {\it sections in $W$ vanishing at
= $P_5, \dots, P_{n+4}$}. The associated projective space then
identifies with the collection of all (general) points in ${\mathcal
M}_1(11)$ giving rise to this point on $S^n$.

\medskip
Put
$$
V_n: = \, {\rm Proj}_{S^n}({\rm coker}({}^tf_n))
$$
where
$$
{}^tf_n: {\mathcal K}_S \boxtimes \dots \boxtimes {\mathcal K}_S
\rightarrow W^\vee \otimes {\mathcal O}_S^\vee.
$$
Then $V_n$ is birational to ${\mathcal M}_1(n+5)$. We have ${\rm
corank}({}^tf_n) \, = \, 6-n$ at a general point of $S^n$. And there
is a natural map
$$\pi: V_n \, \rightarrow \, S^n$$

\medskip

\noindent${\bf n \leq 5}$: \, $\pi$ is surjective with fibres
${\mathbb P}^{5-n}$. Hence $V_n$, which is $\equiv
\overline{\mathcal M}_1(n+5)$, is a {\it rational variety} in this
case. Note that $V_5$ is just $S^5$.

\medskip

\noindent${\bf n=6}$: \, $V = V_6$ is a (reduced) {\it divisor in $S^6$},
hence Gorenstein. It is
defined by the vanishing of det$(f_n)$, which is a section of
${\mathcal K}_{S^n}^{-1}$. So {\bf ${\mathcal K}_V$ is trivial}. We
already know that $h^{(11,0)} =1$ and $h^{(p,0)} = 0$ for $0 < p <
11$ for $\tilde {\mathcal E}^{10}$. These also hold for $V$. So $V$
is Calabi-Yau. Also, the whole construction is rationally defined.

$V=V_6$ fibers over $V_5$ with fibres of dimension $1$, corresponding
to cubics passing through $10$ points. The natural involution on the general fibre,
which is an elliptic curve, gives rise to one, call it $\tau$, on $V$.
The quotient $V/\tau$ is unirational because it is a family of rational curves on a smooth rational surface. Clearly,
$\tau$ must act by $-1$ on the one dimensional space $H^0(V, \Omega^{11})$.

\qed

\bigskip

\section{A C-Y $7$-fold occurring in level $2$}

\medskip

Let $E$ be an elliptic curve with origin $o\in E$, $x\in E$ a point
of order 2 and $y,z\in E$ some other (general) points. Under the
morphism $a:E\to\vert{2[o]+[y]}\vert$, the divisors $2[o]+[y]$ and
$2[x]+[y]$ are linear sections. Of the four points $o$, $x$, $y$ and
$z$ we may assume (under the hypothesis of generality on $y$ and
$z$) that no three are collinear. Thus, we can identify these with
the points $(0:0:1)$, $(1:0:0)$, $(0:1:0)$ and $(1:1:1)$
respectively in order to identify $\vert{2[o]+[y]}\vert$ with
$\mathbb P^2$.

Conversely, let $E$ be a cubic curve in $\mathbb P^2$ which has the
following properties:
\begin{enumerate}
  \item $E$ passes through the points $(0:0:1)$, $(0:1:0)$, $(1:0:0)$
    and $(1:1:1)$.
  \item The line $Y=0$ is tangent to $E$ at the point $(0:0:1)$.
  \item The line $Z=0$ is tangent to $E$ at the point $(0:1:0)$.
\end{enumerate}
Then $E$ is a curve of genus $1$ for which we take $o=(0:0:1)$ as
the origin of a group law. Let $x=(0:1:0)$, $y=(1:0:0)$ and
$z=(1:1:1)$. Then we obtain the relation
\[  2[o] + [y] \simeq 2[x] + [y] \]
It follows that $2x=o$. Thus we have obtained $(E,o,x,y,z)$ of the
type we started with.

Direct calculation shows that the linear system of cubics in
$\mathbb P^2$ that satisfy the conditions above is the linear span
of $X^2Y-XYZ$, $X^2Z-XYZ$, $Y^2Z-XYZ$ and $YZ^2-XYZ$.

\bigskip

Here is an {\bf alternate construction in level $2$}:

\medskip

Let $E$ be an elliptic curve with origin $o\in E$, $x\in E$ a point
of order 2 and $y,z\in E$ some other points. The morphism
$a:E\to\vert{2[o]}\vert$ has fibres $2[o]$, $[y]+[-y]$ and
$[z]+[-z]$ which we map to 0, 1 and $\infty$ respectively in order
to identify $\vert{2[o]}\vert$ with $\mathbb P^1$. The morphism
$b:E\to\vert{[o]+[x]}\vert$ has fibres $[o]+[x]$, $[y]+[x-y]$ and
$[z]+[x-z]$ which we map to 0, 1 and $\infty$ in order to identify
$\vert{[o]+[x]}\vert$ with $\mathbb P^1$. Thus we obtain a morphism
$a\times b:E\to\mathbb P^1\times\mathbb P^1$ which is constructed
canonically from the data $(E,o,x,y,z)$.

Conversely let $E$ be a curve of type $(2,2)$ in $\mathbb
P^1\times\mathbb P^1$ which has the following properties:
\begin{enumerate}
  \item $E$ passes through the points $(0,0)$, $(1,1)$ and
    $(\infty, \infty)$.
  \item The line $\{0\}\times\mathbb P^1$ is tangent to $E$ at the point
    $(0,0)$.
  \item If $(u,0)$ is the residual point of intersection of $E$ with
    $\mathbb P^1\times\{0\}$, then $\{u\}\times\mathbb P^1$ is tangent to $E$ at
    this point.
\end{enumerate}
Then $E$ is a curve of genus 0 for which we take $o=(0,0)$ as the
origin in a group law. Let $x=(u,0)$, $y=(1,1)$ and
$z=(\infty,infty)$. We obtain the identities $2[o]\simeq 2[x]$ from
which it follows that $2x=o$. Thus we have recovered the data
$(E,o,x,y,z)$.

\medskip

\section{Remark on Elliptic curves with Level 3 structure}

\medskip

Let $E$ be an elliptic curve with origin $o\in E$, $x\in E$ a point
of order 3 and $y\in E$ some other point. The morphism
$a:E\to\vert{2[o]}\vert$ has fibres $2[o]$, $[x]+[2x]$ and
$[y]+[-y]$ which we map to 0, 1 and $\infty$ respectively in order
to identify $\vert{2[o]}\vert$ with $\mathbb P^1$. The morphism
$b:E\to\vert{[o]+[x]}\vert$ has fibres $[o]+[x]$, $2[2x]$ and
$[y]+[x-y]$ which we map to 0, 1 and $\infty$ respectively in order
to identify $\vert{[o]+[x]}\vert$ with $\mathbb P^1$. Thus we obtain
a morphism $a\times b:E\to\mathbb P^1\times\mathbb P^1$ which is
constructed canonically from the data $(E,o,x,y)$.

Conversely, let $E$ be a curve of type $(2,2)$ in $\mathbb
P^1\times\mathbb P^1$ which has the following properties:
\begin{enumerate}
  \item $E$ is tangent to the line $\{0\}\times\mathbb P^1$ at the point
    $(0,0)$.
  \item $E$ is tangent to the line $\mathbb P^1\times\{1\}$ at the point
    $(1,1)$.
  \item $E$ passes through the points $(1,0)$ and $\infty,\infty)$.
\end{enumerate}
Then $E$ is a curve of genus $1$ and we use $o=(0,0)$ as the origin
of a group law on $E$. Let $x=(1,0)$, $y=(\infty,\infty)$ and
$p=(1,1)$. We obtain the identities,
\begin{align*}
  2[o] & \simeq [x]+[p] & [o]+[x] & \simeq 2[p]
\end{align*}
It follows that $p=2x$ and $3x=o$. We have thus recovered the data
$(E,o,x,y)$ that we started with.

Direct calculations shows that the linear system of cubic curves in
$\mathbb P^2$ satisfying the conditions given above is the linear
span of $X^2Y-XYZ$, $Y^2Z-XYZ$ and $Z^2X-XYZ$.

\bigskip

\section{Level $3$ and a CY $5$-fold}

\medskip

We will construct a 5-fold with trivial canonical bundle and
singularities only in dimension 2 or less such that its middle
cohomology represents the motive of the (unique) modular form of
level 3 and weight 6.

Consider the linear system $P$ of cubics in ${\mathbb P}^2$ that is
spanned by the curves $X^2Y-XYZ$, $Y^2Z-XYZ$ and $Z^2X-XYZ$; this
system is stable under the cyclic automorphism $X\to Y\to Z\to X$ of
${\mathbb P}^2$. Each curve in the linear system $P$ is tangent to
the line $Z=0$ at the point $p_Y=(0:1:0)$; similarly, the curve is
tangent to $X=0$ at the point $p_Z=(0:0:1)$ and to $Y=0$ at the
point $p_X=(1:0:0)$. Moreover, each curve passes through the point
$p_0=(1:1:1)$. We note that the linear system $P$ is precisely the
collection of cubic curves in ${\mathbb P}^2$ that satisfy these
conditions.

In the divisor class group of a smooth curve in this linear system
we obtain the identities
 \[
    2 p_Y + p_X = 2 p_Z + p_Y = 2 p_X + p_Z = p_X + p_0 + r
 \]
where $r$ denotes the remaining point of intersection of the curve
with the line $Y=Z$ that joins $p_X$ and $p_0$. In particular, we
note that $p_Y-p_X$ is of order $3$ in this class group and $p_Z-p_X
= 2(p_Y-p_X)$.

Conversely, suppose we are given a smooth curve $E$ of genus 1 and a
line bundle $\xi$ of order 3 on $E$; moreover, suppose that three
distinct points $p$, $q$ and $r$ are marked on $E$. We then obtain
two additional points $a$ and $b$ on $E$ such that $a-p=\xi$ and
$b-p=2\xi$ in the divisor class group of $E$. Consider the morphism
$E\to {\mathbb P}^2$ that is given by the linear system of the
divisor $p+q+r$. Moreover, we choose co-ordinates on ${\mathbb P}^2$
so that the point $p$ goes to $p_X$, $q$ goes to $p_0$, $a$ goes to
$p_Y$ and $b$ goes to $p_Z$. This gives an embedding of $E$ as a
curve in ${\mathbb P}^2$ that belongs to the linear system $P$.

Let $S$ denote the surface obtained by blowing up ${\mathbb P}^2$ at
the four points $p_X$, $p_Y$, $p_Z$ and $p_0$, and then further
blowing up the resulting surface at the ``infinitely near points''
that correspond to $Z=0$ at $p_Y$, to $X=0$ at $p_Z$ and to $Y=0$ at
$p_X$. Let $H$ denote the inverse image in $S$ of a general line in
${\mathbb P}^2$; let $E_X$, $E_Y$, $E_Z$ and $E_0$ denote the strict
transforms of the exceptional loci of the first blow-up over the
points $p_X$, $p_Y$, $p_Z$ and $p_0$ respectively; let $F_X$, $F_Y$
and $F_Z$ denote the exceptional divisors of the second blow-up. The
anti-canonical divisor $-K_S=3H - E_X - E_Y - E_Z - E_0 - 2 (F_X -
F_Y - F_Z)$ has a base-point free complete linear system $|-K_S|$
which can be identified with $P$. Let $T$ denote the natural
incidence locus in $S\times P$. The variety
 \[ X = T\times_P T\times_P T \]
is a singular 5-fold which is Gorenstein and has trivial canonical
bundle. Moreover, an open subset of $X_0$ parametrizes tuples of the
form $(E,\xi,p,q,r,s,t,u)$ where $E$ is a curve of genus 1, $\xi$ is
a line bundle of order 3 on $E$ and $p$, $q$, $r$, $s$, $t$ and $u$
are six distinct points on $E$.

Let $L_X$, $L_Y$, $L_Z$ denote the strict transforms in $S$ of the
lines in ${\mathbb P}^2$ defined by $X=0$, $Y=0$, $Z=0$
respectively. Further, let $R$ be the strict transform in $S$ of the
curve in ${\mathbb P}^2$ defined by
 \[ X^2 Z + Y^2 X + Z^2 Y - 3 X Y Z = 0 \]
This is the unique cubic in ${\mathbb P}^2$ that has a node at $p_0$
and is tangent to $X=0$ at $p_Y$, to $Y=0$ at $p_Z$ and to $Z=0$ at
$p_X$. It follows that $R$ is a smooth rational curve that meets
$E_0$ in a pair of distinct points and the triple $(R,L_X,E_Y)$
(respectively $(R,L_Y,E_Z)$ and $(R,L_Z,E_X)$) consists of smooth
curves that meet pairwise transversally.

The morphism $S\to P^*$ induced by the linear system $P$ can be
factorized via a double cover $S\to W$ which is ramified along $R$.
Each of the curves $L_X$ and $E_Y$ (respectively $L_Y$ and $E_Z$;
$L_Z$ and $E_X$) is mapped isomorphically onto the same smooth
irreducible curve $G_Z$ (respectively $G_X$; $G_Y$) in $W$; the
curve $R$ is mapped isomorphically onto the branch locus $Q$ in $W$.
The morphism $W\to P^*$ collapses the curves $G_X$ (respectively
$G_Y$ and $G_Z$) to a point $q_X$ (respectively $q_Y$ and $q_Z$) in
$P^*$; in fact $W\to P^*$ is identified with the blow-up of $P^*$ at
these points. Moreover, $Q$ is mapped to a plane quartic
$\overline{Q}$ which has cusps at these three points.

Let $T$ denote the incidence locus in $P\times S$ as above. It is
the pull-back via $S\to P^*$ of the natural incidence locus
$I\subset P\times P^*$. The latter can be identified (via the
projection $I\to P^*$) with the projective bundle of 1-dimensional
linear subspaces of the tangent bundle of $P^*$. Hence, the
exceptional curve $G_X$ (respectively $G_Y$ and $G_Z$) of the
blow-up $W\to P^*$ can be identified with the fibre $I_X$ of $I\to
P^*$ over the point $q_X$ (respectively $q_Y$ and $q_Z$). Thus we
obtain natural maps $E_{\alpha}\to T$ and $L_{\alpha}\to T$ that are
sections of the ${\mathbb P}^1$-bundle $T\to S$ over the curves
$E_{\alpha}$ and $L_{\alpha}$ respectively; let $\tilde{E_{\alpha}}$
and $\tilde{L_{\alpha}}$ denote the images. Let $T_X$ (respectively
$T_Y$; $T_Z$) denote the fibre of $T$ over the point of intersection
of $L_Y$ and $E_Z$ (respectively $L_Z$ and $E_X$; $L_X$ and $E_Y$).

The tangent direction along $\overline{Q}$ gives a rational morphism
(defined outside the cusps) from $\overline{Q}$ to $I$. It follows
that this extends to a section $R\to T$ of $T\to S$ over $R$ and
gives a curve $\tilde{R}$ in $T$. The quadruple of curves
$\tilde{R}$, $T_X$, $\tilde{L_Y}$, $\tilde{E_Z}$ (respectively,
$\tilde{R}$, $T_Y$, $\tilde{L_Z}$, $\tilde{E_X}$; $\tilde{R}$,
$T_Z$, $\tilde{L_X}$, $\tilde{E_Y}$) meet pairwise transversally in
a single point $r_X$ (respectively $r_Y$; $r_Z$) in $T$. The curve
$\tilde{R}$ in $T$ is mapped to a nodal cubic $\overline{R}$ in $P$
for which $I_X$, $I_Y$ and $I_Z$ are inflectional tangents. The
curves $T_X$, $\tilde{L_Y}$, $\tilde{E_Z}$ (respectively $T_Y$,
$\tilde{L_Z}$, $\tilde{E_X}$; $T_Z$, $\tilde{L_X}$, $\tilde{E_Y}$)
in $T$ lie over $I_X$ (respectively $I_Y$; $I_Z$) in $P$.

The singular locus of the morphism $T\to P$ consists of the curves
$\tilde{R}$, $T_{\alpha}$, $\tilde{L_{\alpha}}$ and
$\tilde{E_{\alpha}}$ for $\alpha=X, Y, Z$ as described above.

The singular fibres of $T\to P$ then have the following description:
\begin{enumerate}
\item If $a$ is a smooth point of $\overline{R}$ which is not a point
  of inflection then the fibre $C_a$ is a rational curve in $S$ with a single
  ordinary node.
\item If $b$ which is on an inflectional tangent (i.~e.\ one one of
  the lines $I_X$, $I_Y$, $I_Z$) of $\overline{R}$ but is {\em not} a
  point of inflection of $\overline{R}$ then the fibre $C_b$ is a
  curve with three components and three nodes (i.~e.\ a ``triangle''
  of ${\mathbb P}^1$'s).
\item If $c$ is a point of inflection of the curve $\overline{R}$,
  then $C_c$ consists of three ${\mathbb P}^1$'s that pass through a point
  and (since $C_c$ lies on a smooth surface $S$) is locally a complete
  intersection.
\item If $d$ is the node of $\overline{R}$ then the fibre $C_d$
  consists of a pair of smooth ${\mathbb P}^1$'s in $S$ that meet in a pair
  of points. In fact the curves are $E_0$ and the strict transform in
  $S$ of the curve in ${\mathbb P}^2$ defined by the equation
  \[    X^2 Y + Y^2 Z + Z^2 X - 3 X Y Z = 0 \]
\end{enumerate}
In particular, the elliptic fibration $T\to P$ is semi-stable but
for the three fibres over the points of inflection of
$\overline{R}$.

Now consider the variety $X=T\times_P T\times_P T$. The singular
points of $X_0$ consist of triples $(x,y,z)$ of points of $T$, where
at least two of these points are critical points for the morphism
$T\to P$. In particular, these points lie over the union of $R$ and
$I_X$, $I_Y$ and $I_Z$. Since the singular points of each of the
fibres described above are isolated, it follows that the singular
locus of $X$ has components of dimension at most 2.

\bigskip

\section{Level $4$ and a C-Y $5$-fold}

\medskip

Let $E$ be an elliptic curve with origin $o\in E$, $x\in E$ a point
of order 4. Under the morphism $a:E\to\vert{3[o]}\vert$, the
divisors $3[o]$, $2[2x]+[o]$, $[2x]+2[x]$ and $[o]+[x]+[3x]$ are
linear sections of the image curve. Let $p$ denote the point of
intersection of the lines corresponding to $3[o]$ and $[2x]+2[x]$.
No three of the points $o$, $2x$, $3x$ and $p$ are collinear. Thus
we can identify $\vert{3[o]}\vert$ with $\mathbb P^2$ in such a way
that $o$ is identified with $(0:1:0)$, $3x$ is identified with
$(0:0:1)$, $p$ is identified with $(1:0:0)$ and $3x$ is identified
with $(1:1:1)$. Thus we obtain a morphism $a:E\to\mathbb P^2$ which
is constructed canonically from the data $(E,o,x)$.

Conversely, let $E$ be a cubic curve in $\mathbb P^2$ which has the
following properties:
\begin{enumerate}
  \item $E$ passes through the points $(0:1:0)$, $(0:0:1)$, $(1:1:1)$
    and $(1:0:1)$.
  \item The line $Z=0$ is an inflectional tangent to $E$ (at the point
    $(0:1:0)$).
  \item The line $X=0$ is tangential to the curve $E$ at the point
    $(0:0:1)$.
  \item the line $Y=0$ is tangential to the curve $E$ at the point
    $(1:0:1)$.
\end{enumerate}
Then $E$ is a curve of genus 0 and we use $o=(0:1:0)$ as the origin
of a group law on $E$. Let $x=(1:0:1)$, $p=(0:0:1)$ and $q=(1:1:1)$.
We obtain the identities
\begin{align*}
  3[0] &\simeq [o]+2[p] &\simeq [p]+2[x] &\simeq [o]+[q]+[x]
\end{align*}
It follows that $2p=o$, $2x=p$ and $q=-x=3x$. Thus we have recovered
the data $(E,o,x)$ that we started with.

Direct calculation shows us that the linear system of cubic curves
in $\mathbb P^2$ that satisfy the above conditions is the linear
span of $YZ(Y-Z)$ and $X(X-Z)^2$.

\medskip

Here is an {\bf alternate construction in level $4$}:

\medskip

Let $E$ be an elliptic curve with origin $o\in E$, $x\in E$ a point
of order 4 and $y\in E$ some other point. The morphism
$a:E\to\vert{2[o]}\vert$ has fibres $2[o]$, $2[2x]$ and $[y]+[-y]$
which we map to 0, 1 and $\infty$ respectively in order to identify
$\vert{2[o]}\vert$ with $\mathbb P^1$. The morphism
$b:E\to\vert{[o]+[x]}\vert$ has fibres $[o]+[x]$, $[2x]+[3x]$ and
$[y]+[x-y]$ which we map to 0, 1 and $\infty$ respectively in order
to identify $\vert{[o]+[x]}\vert$ with $\mathbb P^1$. Thus we obtain
a morphism $a\times b:E\to\mathbb P^1\times\mathbb P^1$ which is
constructed canonically from the data $(E,o,x,y)$.

Conversely, let $E$ be a curve of type $(2,2)$ in $\mathbb
P^1\times\mathbb P^1$ which has the following properties:
\begin{enumerate}
  \item $E$ is tangent to the line $\{0\}\times\mathbb P^1$ at the point
    $(0,0)$.
  \item $E$ is tangent to the line $\{1\}\times\mathbb P^1$ at the point
    $(1,1)$.
  \item If $E$ meets $\mathbb P^1\times\{0\}$ at $(0,0)$ and $(u,0)$ and
  $E$ meets $\mathbb P^1\times\{1\}$ at $(1,1)$ and $(v,1)$; then $u=v$.
  \item $E$ passes through the point $(\infty,\infty)$.
\end{enumerate}
Then $E$ is a curve of genus $0$ and we use $o=(0,0)$ as the origin
of a group law on $E$. Let $x=(u,0)$, $p=(v,1)$ and $q=(1,1)$. We
obtain the identities,
\begin{align*}
  2[o] & \simeq 2[q] &
  [o]+[x] & \simeq  [p]+[q]\\
  2[0] & \simeq  [x]+[p] & \text{~(from condition 3 above)}
\end{align*}
It follows that $q=2x$, $p=3x$ and $4x=o$. Let $y=(\infty,infty)$.
We have thus recovered the data $(E,o,x,y)$ that we started with.
Note that in this construction, $E$ is in $S=
{\mathbb P}^1\times {\mathbb P}^1$, and the Calabi-Yau variety $V$
is a divisor on $S^4=\left({\mathbb P}^1\right)^8$.

\bigskip

\section{Forms of weight $4$ and Calabi-Yau threefolds}

\medskip

The construction of this section is adapted from that of Schoen's
article \cite{Sch}, where he has associated a
modular form $f$ of weight $4$ in this case to a special quintic threefold $X$.
However, we need to
modify his construction to obtain an involution $\tau$ such that the quotient
$X/\tau$ is a rational threefold whose homology is the same as its Chow group.
(Hence, the rational variety has no ``interesting'' motives other than
powers of the Tate motive.) Moreover, the ``additional'' motive on the C-Y
double cover is exactly the one associated (in \cite{Scho}) to the modular form $f$.

Let $E$ be an elliptic curve with $o\in E$ as its origin and $x\in
E$ a point of order 5. Under the morphism $a:E\to\vert{3[o]}\vert$,
the divisors $3[o]$, $[o]+[x]+[4x]$, $2[x]+[3x]$ and $2[3x]+[4x]$
are linear sections. There is a unique identification of
$\vert{3[o]}\vert$ with $\mathbb P^2$ under which these sections are
identified with $Z=0$, $X=0$, $X+Y+Z=0$ and $Y=0$ respectively.

Conversely, let $E$ be  a cubic curve in $\mathbb P^2$ which has the
following properties:
\begin{enumerate}
  \item $E$ passes through the points $(0:1:0)$, $(0:1:-1)$,
    $(1:0:-1)$ and $(0:0:1)$.
  \item The line $Z=0$ is an inflectional tangent to $E$ (at the point
    $(0:1:0)$).
  \item The line $X+Y+Z=0$ is tangent to $E$ at the point $(1:0:-1)$.
  \item The line $Y=0$ is tangent to $E$ at the point $(0:1:-1)$.
\end{enumerate}
The $E$ is a curve of genus 1. Let $o=(0:1:0)$, $x=(0:1:-1)$,
$p=(0:0:1)$ and $q=(1:0:-1)$. We use $o$ as the origin of the group
law on $E$. We obtain the identities,
\[ 3[o] \simeq [o]+[x]+[p] \simeq 2[x]+[q] \simeq 2[q]+[p]. \]
It follows that $q=2x$, $p=4x$ and $5x=o$. Thus we have obtained the
data $(E,o,x)$ that we started with.

Direct calculation shows us that the linear system of cubics in
$\mathbb P^2$ satisfying the above conditions is the linear span of
$YZ(X+Y+Z)$ and $YZ(Y+Z)-X(X+Z)^2$.

Let $E$ be an elliptic curve with $o\in E$ as its origin and $x\in
E$ a point of order 5. The morphism $a:E\to\vert{2[o]}\vert$ has
fibres $2[o]$, $[x]+[4x]$ and $[2x]+[3x]$, which we map to 0, 1 and
$\infty$ to identify $\vert{2[o]}\vert$ with $\mathbb P^1$.
Similarly, the morphism $b:E\to\vert{2[x]}\vert$ has fibres $2[x]$,
$[o]+[2x]$ and $[3x]+[4x]$, which we map to 0, 1 and $\infty$ to
identify $\vert{2[x]}\vert$ with $\mathbb P^1$. Thus we obtain a
morphism $a\times b:E\to\mathbb P^1\times\mathbb P^1$, which is
constructed canonically from the data $(E,o,x)$.

Conversely, let $E$ be a curve of type $(2,2)$ in $\mathbb
P^1\times\mathbb P^1$ which is
\begin{enumerate}
  \item tangent to $\mathbb P^1\times\{0\}$ at the point $(1,0)$,
  \item tangent to $\{0\}\times\mathbb P^1$ at the point $(0,1)$,
  \item passes through the points $(\infty,1)$, $(\infty,\infty)$, and
    $(1,\infty)$.
\end{enumerate}
Then $E$ is a curve of genus 1. Let $o=(0,1)$ and $x=(1,0)$, which
are points on $E$. We use $o$ as the origin of the group law on $E$.
Let $p$, $q$, $r$ denote the points $(\infty,1)$, $(\infty,\infty)$
and $(1,\infty)$ respectively. We obtain the identities,
\begin{align*}
  2[o] & \simeq  [p] + [q] &
  2[x] & \simeq  [q] + [r] \\
  2[o] & \simeq  [x] + [r] &
  2[x] & \simeq  [o] + [p]
\end{align*}
We solve these to show that $a=2x$, $b=3x$, $c=4x$, and $5x=o$. We
have thus recovered the data $(E,o,x)$ that we started with.

Now, appealing to the fibre product paper of Schoen (\cite{Sch}),
we can deduce that our canonical object is Calabi-Yau.

\bigskip

\section{CM forms and C-Y varieties over suitable extensions}

\medskip

Let $E$ be an elliptic curve and $n \geq 1$. Put
$$
B:= \, \{x \in E^{n+1} \, \vert \, \sum_{j=1}^{n+1} x_j = 0\},
$$
which admits an action by the alternating group $A_{n+1}$. Consider
the quotient
$$
X: = \, B/A_{n+1}.
$$
The following result is proved in \cite{PRa}, where the smoothness of
the model for $n=3$ appeals to ideas of Cynk and Hulek.

\medskip

\noindent{\bf Theorem} \, \it $X$ has trivial canonical bundle, with
$H^0(X, \Omega_X^p) = 0$ if $0 < p < n$. If $n \leq 3$, there is a
smooth model $\tilde X$ which is Calabi-Yau. \rm

A submotive $M$ of rank $4$
splits off of $H^3(\tilde X)$ corresponding to sym$^3(H^1(E))$, of
Hodge type
\newline$\{(3,0), (2,1), (1,2), (0,3)\}$. It is simple iff $E$ is
not of CM type, and in this case $\tilde X$ is not rigid.

\medskip

Here is a {\it sketch of proof of Proposition 2}. Let $\Psi$ be the Hecke character
(of weight $k-1$) of an
imaginary quadratic field $K$ attached to $f$, so that $L(s, f) =
L(s, \Psi)$. Pick an algebraic Hecke character $\lambda$ of $K$ of
weight $1$, with an associated elliptic curve $E$ over a finite extension $F_0$.
Then $\Psi/\lambda^{k-1}$ is a finite
order character $\nu$. One attaches (by the Theorem above) a Calabi-Yau variety $V$,
smooth for $k\leq 4$, to $E^{k-1}$. Then it has the requisite properties relative to $f$ over
the number field $F=F_0(\nu)$, having the correct traces of Frobenius elements at the primes outside
a finite set $S$. There are a lot of choices for $\lambda$, which can be used to make
$F$ and $S$ more precise. For example, when the discriminant $-D$ of $K$ is either
odd or divisible by $8$, with $D>4$,
we may choose $\lambda$ to be a {\it canonical} Hecke character of weight one (\cite{Ro}),
which is $K$-valued, equivariant relative to the complex conjugation of $K$,
and is ramified only at the primes dividing $D$. There is an associated CM elliptic curve $E$ over
an explicit number field, which is a quadratic subextension of the Hilbert class field of $K$ when $D$ is odd), studied deeply by Gross
in \cite{Gr}, which is isogenous to all of its Galois conjugates. We start with the C-Y variety $V$ associated
to $E^{k-1}$. In this case $S$ involves only the
primes dividing the level of $f$ and $D$.

\qed

We will discuss Theorem 3 elsewhere
in detail, where for the key case $m=3$ is partly understood via the descent to GSp$(4)/\Q$
(\cite{RaSh}).

\bigskip

\section{Behavior under taking products}

\medskip

In this section we exploit an idea of Claire Voisin. M.~Schuett has remarked
that Cynk and Hulek have also used such an argument.

\medskip

Suppose $X_i$ is a double cover of smooth variety $Y_i$ branched along a
smooth divisor $D_i$ for $i=1,2$; let $\iota_i$ denote the associated
involutions on $X_i$.

The product variety $X_1\times X_2$ carries an action of the
involution $\iota=\iota_1\times \iota_2$ which has $D_1\times D_2$ as its
fixed locus.This is a codimension two transverse intersection of the
divisors $D_1\times X_2$ and $X_1\times D_2$.

Let $X_{12}$ be the blow-up of $X_1\times X_2$ along $D_1\times D_2$ and
$E_{12}$ be the exceptional divisor. Then $E_{12}$ is isomorphic to the
projective bundle over $D_1\times D_2$ associated with the rank two
vector bundle $L_1\oplus L_2$, where $L_i=p_i^* N_{D_i/X_i}$. The strict
transform $E_1$ of $D_1\times X_2$ (respectively $E_2$ of
$X_1\times D_2$) in $X_{12}$ is a divisor that meets $E_{12}$ in a
section of this projective bundle; the two intersections $E_1\cap E_{12}$
and $E_2\cap E_{12}$ are disjoint.

Since $D_1\times D_2$ is the fixed locus of $\iota$, this involution
lifts to $X_{12}$; we denote this lift also by $\iota$ by abuse of
notation. The (scheme-theoretic) fixed locus of $\iota$ on $X_{12}$ is
the smooth divisor $E_{12}$ and hence the quotient $Z_{12}$ of
$X_{12}$ by this involution is a smooth variety. In other words,
$X_{12}$ is a double cover of $Z_{12}$ branched along the smooth divisor
$D_{12}$ which is the image of $E_{12}$.

The involution $\iota_1\times\id_{X_2}$ also lifts to $X_{12}$ since the
base of the blow-up is contained in its fixed locus. Moreover, it
commutes with $\iota$ and hence descends to an involution $\tau$ on
$Z_{12}$. The (scheme-theoretic) fixed locus of this involution on
$Z_{12}$ is the \emph{disjoint} union of the images of $E_1$ and $E_2$
and is thus again a smooth divisor. Thus $Z_{12}$ is the double cover of
a variety $Y_{12}$ branched along a smooth divisor. Moreover, one checks
that $Y_{12}$ is just the blow-up of $Y_1\times Y_2$ along $D_1\times
D_2$ (where by abuse of notation we are using the same notation for the
divisors $D_i$ in $X_i$ and for their images in $Y_i$).

Finally, let us calculate the canonical bundle of $Z_{12}$
and check if it is trivial. First note that the double cover $X$ of
a smooth variety $Y$
branched along a smooth divisor $D$ is obtained by choosing an
isomorphism of $\fO_Y(D)$ with the square of some line bundle $L$ on
$Y$. In this case, the canonical bundle of $X$ is the pull back of
$K_Y\otimes L$.

Next let $L_i$ denote the square root of $\fO_{Y_i}(D_i)$. The canonical
bundle of $K_{X_i}$ is the pull-back of $K_{Y_i}\otimes L_i$; so if
$X_i$ is a Calabi-Yau variety, then $K_{Y_i}$ is the dual of $L_i$.
If $E$ denotes the exceptional locus of the morphism $e:Y_{12}\to
Y_1\times Y_2$, then we see that if $\hat{D_1}$ denotes
the strict transform of $D_1\times Y_2$ in $Y_{12}$ we have
\[ \fO_{Y_{12}}(\hat{D_1}+E) = e^*\fO_{Y_1\times Y_2}(D_1\times Y_2)
\]
Similarly, we have
\[ \fO_{Y_{12}}(\hat{D_2}+E) = e^*\fO_{Y_1\times Y_2}(Y_1\times D_2)
\]
It follows that
\[ \fO_{Y_{12}}(\hat{D_1}+\hat{D_2}) =
    e^* (L_1\otimes L_2\otimes\fO_{Y_{12}}(-E))^2
\]
Since the canonical bundle of $Y_{12}$ is
\[  K_{Y_{12}} = e^*K_{Y_1\times Y_2} \otimes \fO_{Y_{12}}(E)
   =  e^*(L_1\boxtimes L_2)^{-1} \otimes \fO_{Y_{12}}(E) \]
it follows that the canonical bundle of the double cover $Z_{12}$ is
trivial.

\medskip

Theorem 4 now follows by applying this construction.

\qed

\section*{Appendix: \, A consequence of the Hirzebruch Riemann-Roch
theorem}

\medskip

\noindent{\bf Theorem} \, \it The Hirzebruch Riemann-Roch theorem
does not impose any restrictions on the Euler characteristic of an
odd-dimensional smooth and projective variety of dimension at least
3. \rm

\medskip

\noindent{\bf Proof}\, Let $X$ be any smooth projective variety.
Recall that the Todd classes of a variety are the multiplicative
classes defined by the generating function
 \[
    \td(t)=\frac{t}{1-\exp(-t)} =
        1+ \frac{t}{2} + \frac{t^2}{12} -
           \frac{t^4}{720} + \frac{t^6}{30240}
           + O(t^8)
 \]
Fix an integer $m$ and for $i=1,\dots,m$, let $\beta_i$ be algebraic
numbers such that
 \[
    \td(t) \equiv \prod_{i=1}^m (1 + \beta_i t) \mod {t^{m+1}}
 \]
For $X$ of dimension $m$, let $c_i$ for $i=1,\dots,m$ be the Chern
classes of its tangent bundle. Let $\gamma_i$ for $i=1,\dots,m$ be
the Chern roots so that
 \[
    1+c_1 t+c_2 t^2+\dots+c_m t^m = \prod_{i=1}^m (1+\gamma_i t)
 \]
The Todd polynomial of the variety is then given by
 \[
    \Todd(t) = \prod_{i=1}^m td(\gamma_i t)
 \]
The coefficient of $t^m$ in $\Todd(t)$ is the $m$-th Todd class
$\Todd_m$ of the variety. Making use of the above expression
$\td(t)\mod {t^{m+1}}$,
 \[
    \Todd(t) \equiv \prod_{i,j=1}^m (1+\beta_j\gamma_i t)
        \equiv \prod_{j=1}^m \left( 1+ c_1 (\beta_j t) +
                c_2 (\beta_j t)^2 + \dots
                + c_m (\beta_j t)^m \mod{t^{m+1}}
                \right)
 \]
Hence, the coefficient of $c_m$ in the $m$-th Todd class $\Todd_m$
is $\sum_{j=1}^m \beta_j$.

We can compute this as follows. Consider the function
$f(t)=\log(\td(t))$. It has an expression modulo $t^{m+1}$ as
 \[
    f(t) = \sum_{j=1}^m \log(1+\beta_j t) =
        \sum_{j=1}^m \sum_{k=1}^{\infty}
            \frac{(-1)^{k-1}}{k} (\beta_j t)^k
            \mod{t^{m+1}}
 \]
On the other hand we have $f(-t)=f(t)-t$ so that in the expression
of $f(t)$ as a power series in $t$, all the odd degree terms except
$t$ have coefficient 0. In particular, it follows that $\sum_{j=1}^m
\beta_j$ is 0 whenever $m$ is odd and $m>1$.

To summarize, we have proved that the coefficient of the top Chern
class in the Todd class of an odd-dimensional variety is 0. Since
this top Chern class can be identified with the Euler characteristic
we have the result that we wished to prove.

\bigskip

\bibliographystyle{math}    
\bibliography{mf-cy1}

\vskip 0.3in

\end{document}